\documentclass{amsart}
\usepackage{verbatim}
\newtheorem{indh}{Induction Hypothesis}

\newtheorem{pClaim}{{\rm\bf Claim}}
\newcounter{contc}
\newcommand{\enumitoc}{\setcounter{contc}{\value{enumi}}}
\newcommand{\ctoenumi}{\setcounter{enumi}{\value{contc}}}
\newtheorem{theorem}{Theorem}[section]

\newtheorem{lemma}[theorem]{Lemma}

\newtheorem{fact}[theorem]{Fact}

\newtheorem{claim}{Claim}[theorem]
\newtheorem{problem}{Problem}

\newtheorem{Claim}{Claim}

\theoremstyle{definition}

\newtheorem{definition}[theorem]{Definition}

\theoremstyle{remark}

\newcommand{\arablabel}{
          \renewcommand{\labelenumi}{{\rm (\arabic{enumi})}}
          \renewcommand{\theenumi}{{\rm (\arabic{enumi})}}
          \renewcommand{\labelenumii}{{\rm (\arabic{enumii})}}
          \renewcommand{\theenumii}{{\rm (\arabic{enumii})}}
                    }

\newcommand{\alabel}{
          \renewcommand{\labelenumi}{{\rm (\alph{enumi})}}
          \renewcommand{\theenumi}{{\rm (\alph{enumi})}}
          \renewcommand{\labelenumii}{{\rm (\alph{enumii})}}
          \renewcommand{\theenumii}{{\rm (\alph{enumii})}}
                    }

\newcommand{\rlabel}{
          \renewcommand{\labelenumi}{{\rm (\roman{enumi})}}
          \renewcommand{\theenumi}{{\rm (\roman{enumi})}}
          \renewcommand{\labelenumii}{{\rm (\roman{enumii})}}
          \renewcommand{\theenumii}{{\rm (\roman{enumii})}}
                    }
\newcommand{\Rlabel}{
         \renewcommand{\labelenumi}{{\rm(\Roman{enumi})}}
         \renewcommand{\theenumi}{{\rm(\Roman{enumi})}}
         \renewcommand{\labelenumii}{{\rm(\Roman{enumii})}}
         \renewcommand{\theenumii}{{\rm(\Roman{enumii})}}
                    }

\def\myheads#1;#2;{
\pagestyle{myheadings}
\markboth{{\sc\hfill #1\hfill\protect\makebox[0cm][r]{\rm\today}}}
{{\sc\protect\makebox[0cm][l]{\rm\today}\hfill #2\hfill}}
}

\newcommand{\acal}{{\mathcal A}}
\newcommand{\bcal}{{\mathcal B}}
\newcommand{\ccal}{{\mathcal C}}

\newcommand{\fcal}{{\mathcal F}}
\newcommand{\gcal}{{\mathcal G}}
\newcommand{\hcal}{{\mathcal H}}

\newcommand{\ncal}{{\mathcal N}}
\newcommand{\pcal}{{\mathcal P}}

\newcommand{\tcal}{{\mathcal T}}

\newcommand{\ycal}{{\mathcal Y}}

\def\evec{\vec E}

\newcommand{\setm}{\setminus}
\newcommand{\empt}{\emptyset}
\newcommand{\subs}{\subset}

\newcommand{\oo}{{{\omega}_1}}
\newcommand{\dom}{\operatorname{dom}}
\newcommand{\ran}{\operatorname{ran}}
\def\<{\left\langle}
\def\>{\right\rangle}

\def\OO{{\omega}}
\def\oo{\omega_1}
\def\br#1;#2;{\bigl[ {#1} \bigr]^ {#2} }
\def\bc#1;#2;{\bigl( {#1} \bigr)^ {#2} }

\def\ooseq#1;#2;{\< {#1}_{#2}:{#2}<\oo\>}
\def\ooset#1;#2;{\{ {#1}_{#2}:{#2}<\oo\}}
\def\seq#1;#2;#3;{\< {#1}_{#2}:{#2}<#3\>}
\def\set#1;#2;#3;{\{ {#1}_{#2}:{#2}<#3\}}
\def\oseq#1;#2;{\< {#1}_{#2}:{#2}<\OO\>}
\def\oset#1;#2;{\{ {#1}_{#2}:{#2}<\OO\}}
\def\oosequ#1;#2;{\< {#1}^{#2}:{#2}<\oo\>}
\def\oosetu#1;#2;{\{ {#1}^{#2}:{#2}<\oo\}}
\def\sequ#1;#2;#3;{\< {#1}^{#2}:{#2}<#3\>}
\def\setu#1;#2;#3;{\{ {#1}^{#2}:{#2}<#3\}}
\def\osequ#1;#2;{\< {#1}^{#2}:{#2}<\OO\>}
\def\osetu#1;#2;{\{ {#1}^{#2}:{#2}<\OO\}}

\def\force{\raisebox{1.5pt}{\mbox {$\scriptscriptstyle\|$}}
\mbox{$\!\mbox{---}$}}

\newcommand{\fn}{\operatorname{Fn}}

\def\to{\longrightarrow}

\newcommand{\id}{\operatorname{id}}

\def\fin#1;{\br #1;<{\omega};}

\newcommand{\restr}%
{\mathop{\hspace{0.01ex}|\hspace*{-0.02ex}{\grave{}}\hspace{0.4ex}}}

\newcommand{\bijp}{\operatorname{{\rm Bij}_p}}

\newcommand{\isop}{\operatorname{{\rm Iso}_p}}

\newcommand{\adot}{{\dot A}}

\newcommand{\xxdot}{{\dot x}}
\newcommand{\supp}{\operatorname{supp}}
\def\coh#1;{\ccal_{#1}}
\newcommand{\oot}{{\omega}_2}
\def\rnk#1;{\operatorname{rk}_{#1}}
\def\rka{{\rnk {\alpha};}}
\def\rkc{{\rnk {\gamma};}}

\newcommand{\dplus}{\diamondsuit^{+}}

\newcommand{\halfg}{[{\omega};\oo]}
\newcommand{\isom}{\cong}
\def\coh#1;{\ccal_{#1}}

\newcommand{\fdot}{\dot{f}}
\def\enab#1;#2;#3;{$#1$ enables ``$\fdot(#2)=#3$''}

\def\xvec{{\vec x}}
\def\evec{{\vec {\eta}}}
\newcommand{\nest}{\ne^*}
\def\kkvec{{\vec {k}}}

\newcommand{\too}{2^{\oo}}
\newcommand{\ko}{$<\!\!{\omega}$}
\newcommand{\looop}[1]{$#1$-loop}
\newcommand{\ploop}{\looop{p}}
\author{Lajos Soukup}
\address{Mathematical Institute of the Hungarian Academy of Sciences,\\
Budapest, V. Re\'altanoda u. 13-55, Hungary}
\email{soukup@math-inst.hu}
\title{Smooth Graphs}
\thanks{The preparation of this paper was supported by the 
Hungarian National Foundation for Scientific Research grant no. 16391.}
\begin{document}
\begin{abstract}
A  graph $G$ on $\oo$ is called  {\em \ko-smooth }
if for each uncountable $W\subs \oo$, 
$G$ is isomorphic to $G[W\setm W']$ for some finite $W'\subs W$.
 We show that in various models
of ZFC if a graph $G$ is \ko-smooth than G is
necessarily trivial, i.e, either complete or empty. 
On the other hand, we  prove that the existence of a non-trivial,
\ko-smooth graph   is also consistent with ZFC.
\end{abstract}

\maketitle

\section{Introduction}

Answering a question of R. Jamison, H. A. Kierstead and P. J. Nyikos 
proved  in \cite{KN}: {\em 
if the uncountable induced subgraphs of an  uncountable $n$-uniform
hypergraph are pairwise isomorphic,  then the hypergraph  must be either
empty or complete}. 
In this note we investigate how many uncountable subgraphs 
of a graph $G$ on $\oo$ can be isomorphic to $G$ provided that 
it is  {\em non-trivial}, i.e. it 
is not complete or empty. 
As a  corollary of  \cite[theorem 4.2]{HNS} we can
get the following positive result:
the existence  of a non-trivial graph on $\oo$ which embeds into each
of its uncountable subgraphs is consistent with ZFC. 
To formulate  this and the forthcoming results  precisely  
we need the following definition.

\begin{definition}\label{d:smooth}
A  graph $G$ on $\oo$ is called  {\em ${\kappa}$-smooth 
($<\!\!{\kappa}$-smooth)}
if for each uncountable $W\subs \oo$,
$G$ is isomorphic to $G[W\setm W']$ for some $W'\subs W$ 
with $|W'|\le {\kappa}$  ($|W'|<{\kappa})$.
\end{definition}

\begin{fact}
If a graph $G$ on $\oo$ is $n$-smooth for some $n\in {\omega}$,
then $G$  is complete or empty.
\end{fact}

\begin{proof}
Pick ordinals $x_0,x_1,\dots,x_n$ from $\oo$ by finite induction such that 
for each $j\le n$ we have 
\begin{displaymath}
x_{j}\in \bigcap_{i<j}G(x_i)\text{ and }
\biggl|\bigcap_{i\le j}G(x_i)\biggr|=\oo.
\end{displaymath}
If we can not find a suitable $x_j$ then taking
$W=\bigcap_{i<j}G(x_i)$ we have $|W|=\oo$ but $|W\cap G(w)|\le
{\omega}$ for each $w\in W$. Thus $G[W]$ contains an uncountable 
induced empty subgraph and so $G$ is empty.

Assume now that we could choose the sequence $\{x_i:i\le n\}$.
Then let $W=\{x_i:i\le n\}\cup \bigcap_{i\le n}G(x_i)$.
Since $G$ is $n$-smooth there is $W'\subs W$, $|W'|\le n$ such that 
$G\isom G[W\setm W']$. Fix  $i<n+1$ such that $x_i\notin W'$.
Since $x_i\in W\setm W'$, $W\subs G(x_i)\cup \{x_i\}$ and  
$G\isom G[W\setm W']$ 
it follows that there is 
$w\in \oo$ such that $\oo\subs G(w)\cup \{w\}$ and so for each
uncountable $V\subs \oo$ there is $v\in V$ such that
$|V\setm G(v)|\le n$. Thus $G$ contains an uncountable
complete subgraph and so $G$ is complete.
\end{proof}

On the other hand, in \cite[theorem 4.2]{HNS} 
it was shown that $\dplus$ implies that there is a Suslin tree
$\tcal=\<\oo,\prec\>$ such that for each uncountable $X\subs \oo$ there
is a countable $X'\subs X$ such that $\tcal\isom \tcal \restr (X\setm
X')$. Thus the comparability graph  of $\tcal $ is  
${\omega}$-smooth and clearly non-trivial.
However, the question whether  
a \ko-smooth graphs  on $\oo$ is necessarily trivial
was left open. This
gap will be filled up here: we show that (i) 
in different models of ZFC every
\ko-smooth graph on $\oo$ is complete or empty,
(ii) the existence of a non-trivial,
\ko-smooth graph $G$ on $\oo$  is consistent with ZFC. 

The following question however remains unanswered:
\begin{problem}
Is there a non-trivial, ${\omega}$-smooth or
just $\oo$-smooth graph on $\oo$ $($in ZFC$)$? 
\end{problem}

We  use the standard set-theoretical notation throughout, cf
\cite{J}.
For a graph $G$, $V(G)$ denotes the set of vertices of $G$,
$E(G)$ the family of edges of $G$.
If $H\subset V(G)$, $G[H]$  denotes the induced
subgraphs of $G$ on $H$. Given $x\in V(G) $
put $G(x)=\{y\in V(G):\{x,y\}\in E(G)\}$.
If $G$ and $H$ are graphs we write $G\isom H$ 
to mean that $G$ and $H$ are isomorphic. 


If $G$ and $G'$ are graphs,  $\isop(G,G')$ denotes the family 
of isomorphisms between finite induced subgraphs of $G$ and $G'$.

If $q$ is a function let    $\supp(q)=\dom(q)\cup \ran(q)$.

For a cardinal ${\kappa}$ we denote by $\coh {\kappa};$ the standard
poset $\<\fn({\kappa},2;{\omega}),\supseteq\>$
which adds ${\kappa}$ Cohen reals to the ground model.

\section{Models without non-trivial \ko-smooth graphs}
\begin{lemma}
\label{lm:half}
If $G$ is a \ko-smooth graph on $\oo$ and $G$ has 
a --- not necessarily spanned --- subgraph isomorphic to
the bipartite graph $\halfg$ then $G$ is complete.
\end{lemma}

\begin{proof}
Fix $A\in \br \oo;{\omega};$ and $B\in \br \oo;\oo;$
such that $[A,B]\subs E(G)$. 
Let 
\begin{displaymath}
X=\{{\alpha}\in \oo: |\oo\setm G({\alpha})|\le {\omega}\}.    
\end{displaymath}
We show that $X$ is uncountable.
 Indeed, let ${\alpha}<\oo$.
Then for some finite $C\subs A\cup B$ and
$D\subs \oo\setm {\alpha}$ the graphs $G[(A\cup B)\setm C]$ and 
$G[(\oo\setm {\alpha}) \setm D]$ are isomorphic witnessed by a function $f$. 
Then $f''(A\setm C)\subs X$, so $X\not\subs {\alpha}$, i.e.
$|X|=\oo$.

Now, by  recursion, we can  construct a set  
$Y=\{y_{\eta}:{\eta}<\oo\}\cap X$ such that 
$y_{\eta}\in X\cap\bigcap\limits_{{\xi}<{\eta}} G(y_{\xi})$.
Then $G[Y]$ is complete and so $G$ is also complete which was to be proved.
\end{proof}

Let us remark that the statement of  lemma \ref{lm:half} fails for 
${\omega}$-smooth graphs:
the comparability graph $G$ of the Suslin tree $\tcal$ 
constructed in  \cite[theorem 4.2]{HNS}  is  
non-trivial and  ${\omega}$-smooth, but $[{\omega};\oo]\subs G$ and
$[\oo;\oo]\subs \overline{G}$.

Let us recall the definition of splitting number $\mathfrak s$:
\begin{multline*}
\mathfrak s=\min \{|\acal|:{\acal}\subs \br {\omega};{\omega};
\land  \forall X\in \br {\omega};{\omega};\ \exists A\in {\acal}
\ |X\cap A |=|X\setm A|={\omega}\}. 
\end{multline*}

\begin{theorem}\label{tm:no-ko}
Every \ko-smooth graph on $\oo$ is trivial  provided
\ref{kofelt-1} or \ref{kofelt-2} or \ref{kofelt-3} below
hold:
\begin{enumerate}\arablabel
\item \label{kofelt-1} $\oo<\mathfrak s$, 
\item \label{kofelt-2}  $2^{\omega}<2^{\oo}$,
\item \label{kofelt-3}  in a  model  obtained by adding
$\oot$ Cohen reals to some model  $V$.
\end{enumerate}
\end{theorem}

\begin{proof}[Proof of theorem {\rm \ref{tm:no-ko}\ref{kofelt-1}}]
Assume  that $G$ is \ko-smooth.
For each ${\alpha}\in \oo$ let 
$F_{\alpha}=G({\alpha})\cap {\omega}$.
The family $\fcal=\{F_{\alpha}:{\alpha}<\oo\}$ is not a splitting
family for $\mathfrak s>\oo$ so there is an infinite set 
$B\subs {\omega}$ such that $B\subs^* F_{\alpha}$ 
or $B\subs^* {\omega}\setm F_{\alpha}$ for each ${\alpha}\in \oo$.
Then there is $n\in {\omega}$ and an uncountable $I\subs \oo$
such that either $B\setm n\subs F_{\alpha}$ for each ${\alpha}\in I$
or $B\setm n \in {\omega}\setm F_{\alpha}$ for each ${\alpha}\in I$.
Thus either $[B\setm n, I]\subs E(G)$ or
$[B\setm n, I]\cap E(G)=\empt$, i.e $[{\omega};\oo]$
is a subgraph of either $G$ or $\overline{G}$, and so
$G$ is trivial by lemma \ref{lm:half}.   
\end{proof}

\begin{proof}[Proof of theorem {\rm \ref{tm:no-ko}\ref{kofelt-2}}]
Assume on the contrary, that that $G$ is \ko-smooth and non-trivial.
By lemma
\ref{lm:half},  we can choose an uncountable  set $A\subs\oo\setm {\omega}$
such that $ G({\alpha})\cap {\omega}\ne^* G({\beta})\cap {\omega}$
for each $\{{\alpha},{\beta}\}\in \br A;2;$.

For each uncountable $X\subs A$ fix a finite set
$C_X\subs \oo$ and an isomorphism
$f_X$ between $G[({\omega}\cup X)\setm C_X]$ and $G$.
Since $2^{\omega}<2^{\oo}$ there are sets $X,Y\in\br  A;\oo;$
such that $|X\setm Y|\ge{\omega}$, $C_X=C_Y$ and 
$f_X\restr {\omega}= f_Y\restr {\omega}$.
Let ${\xi}\in X\setm Y\setm C_X$.
Then $f= f_Y^{-1}\circ f_X$ is an isomorphism between
$G[({\omega}\cup X)\setm C_X]$ and $G[({\omega}\cup Y)\setm C_Y]$
such that $f\restr ({\omega}\setm C_X)=\id \restr ({\omega}\setm C_X)$. 
Taking ${\eta}=f({\xi})$
we obtain that $G({\xi})\cap ({\omega}\setm C_X)=G({\eta})\cap ({\omega}\setm C_X)$
which contradicts the choice of $A$ because ${\eta}\ne {\xi}$ for
${\xi}\notin \ran(f)$.
\end{proof}

\begin{proof}[Proof of theorem {\rm \ref{tm:no-ko}\ref{kofelt-3}}]
Assume that $G$ is a graph on $\oo$ in $V^{\coh {\oot};}$.
Fix ${\alpha}<\oot$ such that 
$G\in V^{\coh {\alpha};}$. Since 
$\coh {\oot};=\coh {\alpha};*\coh \oot\setm ({\alpha}+\oo);*\coh [{\alpha},{\alpha}+\oo);$, 
by lemma \ref{lm:half} it is enough to prove the following statement:

\begin{lemma}\label{lm:ground-oo}
If $G$ is a graph on $\oo$, $\halfg\not\subs G,\overline{G}$, then $G$ is 
not \ko-smooth in $V^{\coh \oo;}$. 
\end{lemma}

\begin{proof}[Proof of lemma {\rm \ref{lm:ground-oo}}]
 Applying lemma  \ref{lm:half}, we can find an uncountable 
$A\subs  \oo\setm{\omega}$ such that 
$G({\alpha})\cap {\omega}\ne^* G({\beta})\cap {\omega}$ for each 
$\{{\alpha},{\beta}\}\in \br A;2;$.
If $\gcal $ is the $\coh \oo;$-generic filter over $V$,
let $X=\{{\alpha}\in A:\exists p\in \gcal\ p({\alpha})=1\}$.
We show that 
\begin{multline}\notag
\text{$1_{\coh \oo;}\force$ ``{\em $G$
and $G[({\omega}\cup \dot X)\setm Y]$ are not isomorphic}}\\
\text{{\em  for any $Y\in \br \oo;<{\omega};$}. ''}
\end{multline}

Assume on the contrary that 
$p\in \coh \oo;$, $Y\in \br \oo;<{\omega};$ and 
$\fdot$ is a  $\coh \oo;$-name of a function such that
\begin{equation}\notag
\text{
$p\force$ ``$\fdot$ is an isomorphism between 
$G$ and $G[({\omega}\cup \dot X)\setm Y]$''. 
}
\end{equation}
Fix ${\omega}\le {\nu}<\oo$ such that $\dom(p)\cup Y\subs {\nu}$,
$p\force$ ``$\fdot''{\nu}=(({\omega}\cup \dot X)\setm Y)\cap {\nu}$''
and $p\force$ ``$\fdot\restr {\nu}\in V[\gcal\restr {\nu}]$.'' 
From now on we work in $V[\gcal\restr {\nu}]$.
Let $h=f\restr {\nu}$ and  $B=h^{-1}({\omega}\setm Y)$.
Since $G({\alpha})\cap {\omega}\ne^* G({\beta})\cap {\omega}$ for each 
$\{{\alpha},{\beta}\}\in \br A;2;$ it follows that 
if $\{{\zeta},{\xi}\}\in \br \oo\setm {\nu};2;$ then
$G({\zeta})\cap B\ne^*  G({\xi})\cap B$.
Thus for each ${\xi}\in \oo\setm {\nu}$ we have 
$$
f({\xi})={\alpha}\text{ iff }h''(G({\xi})\cap B)=^* (G({\alpha})\cap
{\omega}).
$$ 
Hence $f\restr (\oo\setm {\nu})$ can be defined in $V[\gcal\restr {\nu}]$
and so $X\setm {\nu}\in V[\gcal\restr {\nu}]$, which is impossible by the choice
of $X$.
\end{proof}
The proof of theorem \ref{tm:no-ko} is complete.
\end{proof}

The following theorem claims  that if CH holds in the ground model,
then the   statement of lemma \ref{lm:ground-oo} can be strengthened:
we can find a set  in the  ground model witnessing that $G$ is not 
\ko-smooth in $V^{\coh \oo;}$.

\begin{theorem}\label{tm:oot-miss}
If $CH$ holds and $G$ is a graph on $\oo$
such that $[{\omega};\oo]\not\subs G,\overline{G}$,
then there is an uncountable subset $X$ of $\oo$
such that 
\begin{multline}\notag{}
V^{\coh \oo;}\models
\text{``$G$ is not isomorphic to  $G[X\setm Y]$ 
}
\text{ for any    $Y\in \br \oo;\oo;$''}
\end{multline}
\end{theorem}

The proof is quite long and technical, so we omit it.

\section{Generic construction of a non-trivial \ko-smooth graph}

\begin{theorem}\label{tm:smooth}
If $\too=\oot$, then there is a c.c.c poset $P$ of
size $\oot$ such that 
\begin{displaymath}
V^P\models\text{there is a non-trivial, \ko-smooth graph
$G$ on $\oo$}.
\end{displaymath}
\end{theorem}

\begin{proof}
We construct $P=\ccal*P'$ in two steps: 
in the first step, forcing with $\ccal=\fn(\oo,2;{\omega})$,   
we add  $\oo$-many Cohen reals to $V$ to introduce
our desired graph $G$.
Then, in the second step, we add many isomorphisms between
certain subgraphs of $G$ to $V^\ccal$ to guarantee 
\ko-smoothness of $G$ in $V^{\ccal*P'}$.

To simplify our notation we take $\ccal=\fn(\br \oo;2;,2;{\omega})$
and define the graph $G$ on $\oo$ in $V[\gcal]$, 
where $\gcal$ is the $\ccal$-generic filter over $V$, 
in the straightforward way:
\begin{displaymath}
\text{$\{{\alpha},{\beta}\}\in E(G)$ iff
$\exists p\in \gcal$ $p(\{{\alpha},{\beta}\})=1$.}
\end{displaymath}

If $c\in\ccal$ let $\supp c =\cup\dom  c $ and
$G^c=\<\supp  c ,c^{-1}\{1\}\>$.
Let us remark that if $c,c'\in\ccal$,
$c\le c'$ and $\dom c'=\br \supp c';2;$ then
$G^{c'}$ is a spanned subgraph of $G^c$.

To obtain $P'=P_{\oot}$ we carry out a  finite 
support iteration of c.c.c posets 
$$\<P_{\alpha}:{\alpha}\le\oot ,Q_{\alpha}:{\alpha}<\oot\>$$
in the following way: in the ${\alpha}^{\text{th}}$ step,
we pick an uncountable set $X_{\alpha}$ of $\oo$ in
the intermediate model $V^{\ccal*P_{\alpha}}$ and then we try
to find a finite set $Y_{\alpha}$ and c.c.c poset $Q_{\alpha}$
such that 
\begin{multline}\notag
\text{\em $V^{\ccal*P_{\alpha}*Q_{\alpha}}\models$
``$G$ and $G[X_{\alpha}\setm Y_{\alpha}]$ are isomorphic}\\
\text{\em witnessed by a function $f_{\alpha}$.''}
\end{multline}

The poset $Q_{\alpha}$ will consist of certain 
isomorphisms between finite subgraphs of $G$ and
$G[X_{\alpha}\setm Y_{\alpha}]$, ordered by the 
reverse inclusion.
In other words, we force with certain finite approximations
of an isomorphism between 
$G$ and $G[X_{\alpha}\setm Y_{\alpha}]$. 

The problem is the right choice of $Q_{\alpha}$ because
we should meet two contradictory requirements. First, 
the poset $Q_{\alpha}$ should satisfy c.c.c
and forcing with $Q_{\alpha}$ can not introduce an uncountable empty
or complete subgraph  of $G$, 
therefore $Q_{\alpha}$ can not contain
too many elements. 
On the other hand, to guarantee that a $Q_{\alpha}$-generic filter gives an isomorphism
between $G$ and $G[X_{\alpha}\setm Y_{\alpha}]$
we need some density arguments, i.e. certain 
subsets of $Q_{\alpha}$ should be dense in $Q_{\alpha}$, which involves
that $Q_{\alpha}$ can not be too small.
As it turns out, it will be quite easy
to meet the first requirement, the hard part of the proof is 
how to cope with the second one.

Now assume that $P_{\alpha}$ is constructed and 
let us see the induction step.

First, using a bookkeeping function, we  pick the set 
$X_{\alpha}\in \br \oo;\oo;\cap V^{\ccal*P_{\alpha}}$ in such a way  that
\begin{equation}
\tag{$*$}
\{X_{\alpha}:{\alpha}<\oot\}=\br \oo;\oo;\cap V^{\ccal*P_{\oot}}.
\end{equation}

To construct the poset $Q_{\alpha}$ we need the following
induction hypothesize.
To formulate it we use two  notions.
A graph  $G$  
is {\em strongly non-trivial} provided that each uncountable 
family  of pairwise
disjoint, finite subsets of $V(G)$  
contains four distinct elements,
a, b, c, d such that  $[a,b]\subset E(G)$ and 
$[c,d]\cap E(G) =\empt$.
If $G$ is a graph,  
a set $A\subs V(G)$ is called {\em dense in G\ } iff for each
pair $B$ and $B'$ of disjoint finite subsets of  $V(G)$
there is  ${\alpha}\in A$ such that
$G({\alpha})\supset B$ and $ G({\alpha})\cap B'=\empt$.

\begin{indh}\makebox[1pt]{}
\begin{enumerate}\Rlabel
\item \label{non-trivial}
$V^{\ccal*P_{\alpha}}\models$ ``$G$ is strongly non-trivial'',
\item \label{many-dense} $V^{\ccal*P_{\alpha}}\models$ ``
$\forall X\in\br \oo;\oo;$ 
$\exists Y\in \br X;<{\omega};$ $\forall {\delta}<\oo$ 
$\exists A\in \br X\setm {\delta};{\omega};$ $A$ is dense in 
$G[X\setm Y]$'',
\end{enumerate}
\end{indh}

The preservation of the induction hypothesize \ref{non-trivial} and
\ref{many-dense} during the iteration 
will be verified  later 
in lemmas \ref{lm:non-trivial} and
\ref{lm:many-dense}.

We continue the construction of the poset $Q_{\alpha}$.
Using \ref{many-dense} fix $Y_{\alpha}\in \br X_{\alpha};<{\omega};$ 
and pairwise disjoint countable subsets 
$\{D_{\xi}:{\xi}<\oo\}$ of $X_{\alpha}\setm Y_{\alpha}$
which are dense in $G[X_{\alpha}\setm Y_{\alpha}]$.

Let us recall that for each ${\beta}<{\alpha}$ 
in the ${\beta}^{\text{th}}$ step
we already constructed 
an  isomorphism $f_{\beta}$ between $G$
and $G[X_{\beta}\setm Y_{\beta}]$.
For each ${\beta}<{\alpha}$ the set 
$C_{\beta}=\{{\nu}<\oo:f_{\beta}{}''{\nu}\subs {\nu}\}$ is
clearly  club and 
$C_{\beta}$ belongs to $V^{\ccal*P_{{\beta}}*Q_{\beta}}\subs V^{\ccal*P_{\alpha}}$. 
Since $P_{\alpha}$ satisfies c.c.c and 
$|{\alpha}|<\too={\omega}_2$,  there is a club set
$C\subs \oo$ even in $V$ such that $|C\setm C_{\beta}|\le {\omega}$ for each
${\beta}<{\alpha}$. 

The club set $C=\{{\gamma}_{\nu}:{\nu}<\oo\}$ 
gives a natural partition 
$\acal_{\alpha}=\{A^{\alpha}_{\nu}:{\nu}<\oo\}$ of $\oo$
into countable pieces: let
$A^{\alpha}_{\nu}=[{\gamma}_{\nu},{\gamma}_{{\nu}+1})$
for ${\nu}<\oo$. 
We can thin out  $C$ to contain
only limit ordinals and in this case every $A^{\alpha}_{\nu}$
is infinite.
Define the map $\rka:\oo\to \oo$ by the formula 
${\xi}\in A^{\alpha}_{\rka(\xi)}$.

If ${\beta}<{\alpha}$ then $|C\setm C_{\beta}|\le {\omega}$ 
and so  all but countably many 
$A^{\alpha}_{\nu}$'s are $f_{\beta}$-closed.
By shrinking $C$ we can assume every $A^{\alpha}_{\eta}$ 
contains some $D_{\xi}$
and so 
\begin{enumerate}\rlabel
\item \label{xdense}
$A^{\alpha}_{\eta}\cap (X_{\alpha}\setm Y_{\alpha})$
is dense in $G[X_{\alpha}\setm Y_{\alpha}]$.
\enumitoc
\end{enumerate}
Since  $A^{\alpha}_{\eta}\in V$ and infinite, it follows  
\begin{enumerate}\rlabel
\ctoenumi
\item \label{odense}
$A^{\alpha}_{\eta}$ is dense in $G$.
\enumitoc
\end{enumerate}

For ${\eta}<\oo$ let $O_{\eta}=[{\omega}{\eta},{\omega}{\eta}+{\omega})$ and 
 $B^{\alpha}_{\eta}=\bigcup\{A^{\alpha}_{\eta}:{\nu}\in O_{\eta}\}$.
Put $\bcal_{\alpha}=\<B^{\alpha}_{\eta}:{\eta}<\oo\>$.

Given two sets $Z$ and $W$ denote by $\bijp(Z,W)$ the family of
bijections between finite subsets $Z$ and $W$.

If $p\in \bijp(\oo,X\setm Y)$  
a sequence 
$
\xvec=\<x_0,x_1,\dots,x_n,\>
$ 
of countable ordinals is a {\em \ploop}  iff   
$n\ge 1$, $x_0=x_n$ and there is a sequence 
$\<k_0,\dots,k_{n-1}\>\in {}^n \{-1,+1\}$
such that 
\begin{enumerate}\rlabel\ctoenumi
\item \label{p-rank}
$\rnk {\alpha};(x_{i+1})=\rnk {\alpha}; (p^{k_i}(x_i))$ for each
$i<n$, 
\item \label{p-triv}
there is no $i<n$ such that $\{k_i,k_{i+1}\}=\{-1,+1\}$,
$x_{i+1}=p^{k_i}(x_i)$ and $x_{i+2}=p^{k_{i+1}}(x_{i+1})$.
\enumitoc
\end{enumerate}

We say that  $p$ is {\em loop-free} if there is no 
$p$-loop.

Now we are in the position to define the poset $Q_{\alpha}$.
We put a  finite function 
$p\in\isop(G,G[X_{\alpha}\setm Y_{\alpha}])$ into $Q_{\alpha}$
iff 
\begin{enumerate}\rlabel
\ctoenumi
\item \label{p-closed}
$p'' B_{\eta}\subs B_{\eta}$ for each ${\eta}<\oo$,
\item \label{p-loop-free}
$p$ is loop-free.
\enumitoc
\end{enumerate}
As promised, $Q_{\alpha}$ is ordered by the reverse inclusion:
$Q_{\alpha}=\<Q_{\alpha},\supseteq\>$.

Let us recall that $\supp p=\dom(p)\cup\ran(p)$
for $p\in Q_{\alpha}$.

We need to show that $Q_{\alpha}$ satisfies c.c.c and a
$Q_{\alpha}$-generic filter gives an isomorphism 
between $G$ and $G[X_{\alpha}\setm Y_{\alpha}]$.
First we prove an auxiliary lemma.

\begin{lemma}\label{lm:disj_un}
If $p,q\in \bijp(\oo,\oo)$, $\rka''\supp p\cap \rka''\supp q=\empt$
and $\xvec=\<x_0,\dots,x_n\>$ is a \looop{(p\cup q)}, 
then $\xvec$ is either a $p$-loop
or a $q$-loop.
\end{lemma}
\begin{proof}
Assume that $x_0\in\supp p$. Then 
$x_0\notin\supp q$, so $\rka (x_1)=\rka (p^{k_0}(x_0))$ for some
$k_0\in \{-1,+1\}$. Since $p^{k_0}(x_0)\in \supp p$ we have 
$\rka (x_1)=\rka (p^{k_0}(x_0))\notin\rka''\supp q$ and so
$x_1\notin\supp q$. Repeating this argument  we yield
$\{x_0,\dots,x_n\}\subs \supp p\setm \supp q$ and so $\xvec$ is a
\looop{p}. 
\end{proof}

\begin{lemma}\label{lm:ccc}
$Q_{\alpha}$ satisfies c.c.c.
\end{lemma}

\begin{proof}
We work in $V^{\ccal*P_{\alpha}}$.
Assume that 
 $\{q_{\xi}:{\xi}<\oo\}\subs Q_{\alpha}$,
$c_{\xi}=\supp q_{\xi}$ and $r_{\xi}=\rka''c_{\xi}$. 
Applying standard $\Delta$-system and counting arguments 
we can find $I\in \br \oo;\oo;$ such that 
\begin{enumerate}\arablabel
\item $\{c_{\xi}:{\xi}\in I\}$ forms a $\Delta$-system with kernel
$c$,
\item $\{r_{\xi}:{\xi}\in I\}$ forms a $\Delta$-system with 
kernel $r$,
\item   $\rka''c=r$,
\item  $\rka''(c_{\xi}\setm c)=r_{\xi}\setm r$ for each ${\xi}\in I$,
\item  $q_{\xi}\restr c=q'$ for each ${\xi}\in I$.  
\end{enumerate}
 

Since $G$ is strongly non-trivial in $V^{\ccal*P_{\alpha}}$
by the induction hypothesis \ref{non-trivial}, 
there is $\{ {\xi},{\zeta}\}\in \br I;2;$ 
such that $[c_{{\xi}}\setm c,c_{{\zeta}}\setm c]\subs E(G)$.
We show that $q=q_{{\xi}}\cup q_{{\zeta}}\in  Q_{\alpha}$.
Clearly $q\in\isop(G, G[X_{\alpha}\setm Y_{\alpha}])$ and 
$q$  satisfies
\ref{p-closed}. Since $q=q'\cup (q_{{\xi}}\setm q')\cup 
(q_{{\zeta}}\setm q') $ and 
the sets $\rka'' q'$, $\rka''(q_{\xi}\setm q')$ and 
$\rka'' (q_{\zeta}\setm q')$ are 
pairwise disjoint we have that $q$ satisfies \ref{p-loop-free} as well
by lemma \ref{lm:disj_un}.
\end{proof}

If $\gcal^{Q_{\alpha}}$ is the $Q_{\alpha}$-generic filter over 
$V^{\ccal*P_{\alpha}}$ let
$f_{\alpha}=\cup\{q:q\in\gcal^{Q_{\alpha}}\}$.

\begin{lemma}\label{lm:dr-dense}
$V^{\ccal*P_{\alpha}*Q_{\alpha}}\models$
``$f_{\alpha}$ is an isomorphism between $G$ and 
$G[X_{\alpha}\setm Y_{\alpha}]$.''
\end{lemma}

\begin{proof}
We need to prove that $\dom(f_{\alpha})=\oo$ 
and $\ran(f_{\alpha})=X_{\alpha}\setm Y_{\alpha}$ which 
follows if for each  
${\nu}\in \oo$ and ${\mu}\in X\setm Y$ both 
$$
D_{\nu} =\{q\in Q_{\alpha}:{\nu}\in\dom q\}
$$
and
$$
R_{\mu} =\{q\in Q_{\alpha}:{\mu}\in \ran q\}
$$
are dense in $Q_{\alpha}$.
Fix $q\in Q_{\alpha}$.
Write $\rka({\nu})={\omega}{\eta}+n$.
Pick  ${\omega}{\eta}\le {\zeta}<{\omega}{\eta}+{\omega}$ such that 
$(\supp q)\cap A^{\alpha}_{\zeta}=\empt$.
Since $A^{\alpha}_{\zeta}\cap (X_{\alpha}\setm Y_{\alpha})$ is dense in 
$G[X_{\alpha}\setm Y_{\alpha}]$ we can find
${\nu}' \in A^{\alpha}_{{\zeta}}\cap (X_{\alpha}\setm Y_{\alpha})$ such that
$\{{\nu}',q({\xi}) \}\in E(G)$ iff $\{{\nu},{\xi}\}\in E(G)$ 
for each ${\xi}\in\dom q$. Let 
$q'=q\cup\{\<{\nu},{\nu}'\>\}$. By the choice of ${\zeta}'$,
$\rka({\nu}')={\zeta}\notin \rka''(\supp q)$, so this extension
of $q$ can not introduce a $q'$-loop, i.e. $q'\in Q_{\alpha}$.
Thus $q'\in D_{{\nu}}$ and  $q'\le q$ which was to be proved.
The density of $R_{\mu}$ can be verified by a similar argument
using the density of $A^{\alpha}_{\zeta}$ in $G$. 
\end{proof}

The induction step is complete so the theorem is proved 
provided we can verify the induction hypothesize 
 \ref{non-trivial} and
\ref{many-dense}  in every $V^{\ccal*P_{\gamma}}$.
First we deal with \ref{non-trivial} because it is fairly easy.
Checking \ref{many-dense} is the crux of our proof.

\begin{lemma}\label{lm:non-trivial}
The induction hypothesis \ref{non-trivial} holds, i.e. 
$G$ is strongly non-trivial in every 
$V^{\ccal*P_{{\alpha}}}$.
\end{lemma}

\begin{proof}
First remark that $G$ is clearly strongly non-trivial in 
$V^{\ccal}$. 
By  \cite[lemma 4.10]{HNS} we can assume that 
${\alpha}={\gamma}+1$ and $G$ is strongly non-trivial
in $V^{\ccal*P_{\gamma}}$.
Working in $V^{\ccal*P_{\alpha}}$ assume that 
$q\force $ ``{\em $\{\xxdot_{\xi}:{\xi}<\oo\}$ are pairwise
disjoint, finite subsets of $\oo$.}''
For each ${\xi}<\oo$ pick a condition $q_{\xi}\le q$ and 
a finite subset $x_{\xi}$ of $\oo$
such that $q_{\xi}\force$ ``$\xxdot_{\xi}=x_{\xi}$''.
Since $Q_{\gamma}$ satisfies c.c.c, we can assume that the
sets $x_{\xi}$ are pairwise disjoint.

We can also assume that $x_{\xi}\subs \dom q_{\xi}$ because in 
lemma \ref{lm:dr-dense} we showed that the sets 
$D_{\nu}$ are dense in $Q_{\gamma}$.
 
From now on we can argue as in lemma \ref{lm:ccc}.
Let 
$c_{\xi}=\supp q_{\xi}$ and $r_{\xi}=\rkc''c_{\xi}$. 
We can find $I\in \br \oo;\oo;$
such that $\{c_{\xi}:{\xi}\in\}$ forms a $\Delta$-system with kernel
$c$ and $\{r_{\xi}:{\xi}\in I\}$ forms a $\Delta$-system with 
kernel $r$, moreover 
$\rkc''c=r$, $\rkc''(c_{\xi}\setm c)=r_{\xi}\setm r$,  
$q_{\xi}\restr c$ is independent from ${\xi}$  and 
$x_{\xi}\subs c_{\xi}\setm c$ for 
each ${\xi}\in I$. Write $c'_{\xi}=c_{\xi}\setm c$,
$q'_{\xi}=q_{\xi}\restr c'_{\xi}$, $r'_{\xi}=r_{\xi}\setm r$
and $q'=q_{\xi}\restr c$.

Since $G$ is strongly non-trivial in $V^{\ccal*P_{\xi}}$ there
are ${\xi}_0, {\xi}_1,{\zeta}_0, {\zeta}_1\in I$ 
such that $[c'_{{\xi}_0},c'_{{\zeta}_0}]\subs E(G)$ and
$[c'_{{\xi}_1},c'_{{\zeta}_1}]\cap E(G)=\empt$ .
Then $q^i=q_{{\xi}_i}\cup q_{{\zeta}_i}\in \isop(G, G[X\setm Y])$ and 
$q^i$ clearly satisfies
\ref{p-closed}. Since $q^i=q'\cup q'_{{\xi}_i}\cup q'_{{\zeta}_i} $ and 
the sets $\rkc''q'$, $\rkc''q'_{\xi_i}$ and $\rkc'' q'_{\zeta_i}$ are 
pairwise disjoint we have that $q^i$ satisfies \ref{p-loop-free} as well by
lemma \ref{lm:disj_un}.
Thus
$$q^0\force [\xxdot_{{\xi}_0},\xxdot_{{\zeta}_0}]\subs E(G)$$
and
$$q^1\force [\xxdot_{{\xi}_1},\xxdot_{{\zeta}_1}]\cap E(G)=\empt.$$
\end{proof}

Now we start to work on (II).
\begin{definition}\label{df:act-freely}
Assume that $\hcal$ is a family of function, 
$\dom(h)\cup\ran(h)\subs \oo$ for each $h\in\hcal$. 
A sequence $\xvec=\<x_0,x_1,\dots, x_n\>\in {}^n\oo$ is called 
{\em $\hcal$-loop} if $n\ge 1$, $x_0=x_n$, and there are  sequences
$\<h_0,\dots,h_{n-1}\>\in {}^n\hcal$ and 
$\<k_0,\dots,k_{n-1}\>\in {}^n \{-1,+1\}$ such that  
\begin{enumerate}\rlabel
\ctoenumi
\item \label{h-rank}
$h_i^{k_i}(x_i)=x_{i+1}$ for each $i<n$,
\item \label{h-triv}
there is no $i<n-1$ such that $h_i=h_{i+1}$ and
$\{k_i,k_{i+1}\}=\{-1,+1\}$.
\enumitoc
\end{enumerate}
Let $Z\subs \oo$. We say that {\em $\hcal$ acts loop-free on $Z$}
if
\begin{enumerate}\rlabel
\ctoenumi
\item $Z$ is $h$-closed for each $h\in\hcal$,
\item $Z$ does not contain any $\hcal$-loop.
\enumitoc
\end{enumerate}
\end{definition}
\begin{definition}\label{df:det}
A condition $p=\<c,q\>\in \ccal*P_{\alpha}$ is called 
{\em determined} iff 
\begin{enumerate}\arablabel
\item $q$ is a function, $\dom(q)\in \br \oo;<{\omega};$, 
\item $q({\eta})$ is a function for each ${\eta}\in\dom(q)$,
\item $\bigcup\{\supp q({\eta}):{\eta}\in \dom (q)\}\subs \supp c$,
\item $\dom(c)=\br \supp c;2;$.
\end{enumerate}
\end{definition}
The determined conditions are dense in $\ccal*P_{\alpha}$.

\begin{lemma}\label{lm:loop-free}
In $V^{\ccal*P_{{\alpha}}}$
for each  $J\in \br {\alpha};<{\omega};$ there is ${\mu}<\oo$
such that $\{f_{\xi}:{\xi}\in J\}$ acts loop-free on $\oo\setm {\mu}$.
\end{lemma}

\begin{proof}

We work in $V[\gcal]$, where $\gcal$ is the $\ccal*P_{\alpha}$-generic
filter over $V$.
The lemma will be proved by induction on $\max J$. 
Let ${\zeta}=\max J$ and $J'=J\setm \{{\zeta}\}$.
Using the inductive hypothesis fix ${\mu}<\oo$ such that 
\begin{enumerate}\alabel
\item ${\mu}=\bigcup\{B\in \bcal_{\zeta}:B\cap{\mu}\ne \empt\}$,
\item if $A\in\acal_{\zeta}$ and $A\subs \oo\setm {\mu}$ then
$A$ is $f_{\xi}$-closed for each ${\xi}\in J'$,
\item $\{f_{\xi}:{\xi}\in J'\}$ acts loop-free on $\oo\setm {\mu}$.
\end{enumerate}
Assume on the contrary that  
$\<x_0,\dots,x_n\>\in {}^n(\oo\setm {\mu})$ is an 
$\{f_{\xi}:{\xi}\in J\}$-loop witnessed by the sequences
$\<g_i:i<n\>\in {}^n\{f_{\xi}:{\xi}\in J\}$ and 
$\<k_i:i<n\>\in {}^n \{-1,+1\}$.
Let $M=\{m<n:g_m=f_{\zeta}\}$. By the induction hypothesis
$M\ne\empt$.
Write $M=\{m_j:j<\ell\}$, $m_0<\dots<m_{\ell-1}$.
Let $y_0=x_{m_0}$, $y_1=x_{m_1}$, $\dots$, $y_{\ell-1}=x_{m_{\ell-1}}$
and $y_{\ell}=x_{m_0}$.
Pick a determined condition $\<c,q\>\in \gcal$ such that 
 $y_j,f_{\zeta}^{k_{m_j}}(y_j)\in
\dom(q({\zeta}))\cap\ran(q({\zeta}))$ for each $j<\ell$. 
We claim that $\<y_j:j\le \ell\>$ is a $q({\zeta})$-loop 
witnessed by the sequence $\<k_{m_j}:j<\ell\>$, which contradicts the choice
of $Q_{\zeta}$. 
Condition \ref{p-rank} holds because  
$\rnk {\zeta};(y_{j+1})=\rnk {\zeta};(f_{\zeta}^{k_{m_j}}(y_j))$ by (b).
Assume on the contrary that \ref{p-triv} fails, i.e, there is $j<\ell$
such that $\{k_{m_j},k_{m_{j+1}}\}=\{-1,+1\}$,
$y_{j+1}=f^{k_{m_{j}}}_{\zeta}(y_j)$ and  
$y_{j+2}=f^{k_{m_{j+1}}}_{\zeta}(y_{j+1})$.
Since $f^{k_{m_{j}}}_{\zeta}(y_j)=f^{k_{m_{j}}}_{\zeta}(x_{m_j})=
x_{m_j+1}$ and $y_{j+1}=x_{m_{j+1}}$, and so $x_{m_j+1}=x_{m_{j+1}}$,
 by (c) it follows that 
$m_{j}+1=m_{j+1}$. Similarly, $m_{j+1}+1=m_{j+2}$.
Thus $x_{m_j}=y_j$, $x_{m_j+1}=y_{j+1}$ and $x_{m_j+2}=y_{j+2}$.
So $g_{m_j}=g_{m_j+1}=f_{\zeta}$ and 
$\{k_{m_j},k_{m_j+1}\}=\{-1,+1\}$ which contradicts our assumption
that $\<g_i:i<n\>$ and $\<k_i:i<n\>$ satisfied \ref{h-triv}.
\end{proof}

\begin{lemma}\label{lm:many-dense}
The induction hypothesis \ref{many-dense} holds
in $V^{\ccal*P_{\alpha}}$, i.e.\\
$V^{\ccal*P_{\alpha}}\models$ ``
$\forall X\in\br \oo;\oo;$ 
$\exists Y\in \br X;<{\omega};$ $\forall {\delta}<\oo$ 
$\exists A\in \br X\setm {\delta};{\omega};$ $A$ is dense in 
$G[X\setm Y]$'',
\end{lemma}

\begin{proof}

Assume that
$$
1_{\ccal*P_{{\alpha}}}\force X=\{\xxdot_{\xi}:{\xi}<\oo\}\in \br \oo;\oo;.
$$
Pick determined conditions $p_{\xi}=\<c_{\xi},q_{\xi}\>\in \ccal*P_{{\alpha}}$
and $x_{\xi}\in\oo$
such that $p_{\xi}\force$ ``$\xxdot_{\xi}=x_{\xi}$''. 
We can assume that $x_{\xi}\in \supp c_{\xi}$. 
Write  $J_{\xi}=\dom q_{\xi}$ and $Z_{\xi}=\supp(c_{\xi})$.

Now there is $K\in \br \oo;\oo;$ such that the conditions
$\{p_{\xi}:{\xi}\in K\}$ are ``{\em pairwise twins}'', i.e.  
\begin{enumerate}\arablabel
\item $\{Z_{\xi}:{\xi}\in K\}$ forms a $\Delta$-system with kernel $Z$,
\item  $\{J_{\xi}:{\xi}\in K\}$ forms a $\Delta$-system with kernel $J$,
\item $\max Z <\min (Z_{\xi}\setm Z) < \max (Z_{\xi}\setm Z) <\min
(Z_{\xi'}\setm Z)$
for ${\xi}<{\xi}'\in K$,
\item $|Z_{\xi}|=|Z_{\xi'}|$ for $\{{\xi},{\xi}'\}\in \br K;2;$.
Denote by ${\varphi}_{{\xi},{\xi}'}$ the natural bijection between
$Z_{\xi}$ and $Z_{\xi'}$.
\item $c_{\xi'}
(
(\{{\varphi}_{{\xi},{\xi'}}({\nu}),
{\varphi}_{{\xi},{\xi'}}({\nu}')\}))=
c_{\xi}(\{{\nu},{\nu}'\})$
for $\{{\nu},{\nu}',\}\in \br Z_{\xi};2;$ and 
$\{{\xi},{\xi}'\}\in \br K;2;$,
\item $q_{\xi'}({\eta})=\{\<{\varphi}_{{\xi},{\xi'}}({\nu}),
{\varphi}_{{\xi},{\xi'}}({\nu'})\>:\<{\nu},{\nu}'\>\in
q_{{\xi}}({\eta})\}$ for ${\eta}\in J$ and 
$\{{\xi},{\xi}'\}\in \br K;2;$.
\end{enumerate}

Since  $\bcal_{\eta}$ is a partition of $\oo$ into countable 
pieces for ${\eta}\in J$,
 there is a club set 
$C=\{{\gamma}_{\nu}:{\nu}<\oo\}\subs \oo$ in 
$V^{\ccal*\pcal_{{\alpha}}}$
such that for each  ${\eta} \in J$ and ${\nu}<\oo$ we have
$$
[{\gamma}_{\nu},{\gamma}_{{\nu}+1})=\bigcup\{B\in\bcal_{\eta}:
B\cap [{\gamma}_{\nu},{\gamma}_{{\nu}+1})\ne\empt\}.   
$$ 

Since $\ccal*P_{\alpha}$ is c.c.c we can assume that $C\in V$.

By thinning out $K$ we can assume that if ${\xi}<{\xi}'\in K$ then
there is ${\gamma} \in C$ such that 
$\max (Z_{\xi}\setm Z)<{\gamma}<\min (Z_{{\xi}'}\setm Z)$, moreover
$\max Z<\min C$. 

By lemma \ref{lm:loop-free} fix ${\mu}\in C$ such that ${\delta}\le{\mu}$ and 
$1_{\ccal*P_{{\alpha}}}\force$  ``{\em $\{f_{\eta}:{\eta}\in J\}$ 
acts loop-free on $\oo\setm {\mu}$}''. 

If  $\evec=\<{\eta}_0,\dots, {\eta}_{n-1}\>\in {}^nJ$ and 
$\kkvec=\<k_0,\dots,k_{n-1}\>\in {}^n\{-1,+1\}$ for some 
$n\in {\omega}$ then
let
\begin{displaymath}
f_{\<\evec,\kkvec\>}=f_{{\eta}_{n-1}}^{k_{n-1}}\circ
\dots\circ f_{{\eta}_{0}}^{k_{0}}.
\end{displaymath}

If $p=\<c,q\>$ is determined and $J\subs \dom (q)$ 
we define the {\em $q$-approximation of $f_{\<\evec,\kkvec\>}$},
$f^{q}_{\<\evec,\kkvec\>}$, in the natural way:
\begin{displaymath}
f^{q}_{\<\evec,\kkvec\>}=q({\eta}_{n-1})^{k_{n-1}}\circ
\dots\circ q({\eta}_{0})^{k_{0}}.
\end{displaymath}

We say that $f_{\<\evec,\kkvec\>}$ is {\em irreducible} if there is no 
$i<n-1$ such that ${\eta}_i={\eta}_{i+1}$ and $\{k_i,k_{i+1}\}=\{-1,+1\}$.

Let ${\xi}\in K$ be arbitrary.
An irreducible $f_{\<\evec,\kkvec\>}$ is {\em active} iff 
$\dom f^{q_{\xi}}_{\<\evec,\kkvec\>}\cap (Z_{\xi}\setm Z)\ne\empt$, i.e., 
there is a sequence 
$\xvec=\<x_0,\dots, x_{n-1}\>\in {}^n (Z_{\xi}\setm Z)$ such that 
 $x_{i+1}=q_{{\xi}}({\eta}_i)^{k_i}(x_i)$ for $i<n$.
Observe that the definition of activeness above does not depend on
the choice ${\xi}$
because the conditions $\{\<c_{\xi},q_{\xi}\>:{\xi}\in K\}$
are pairwise twins.

We say that $\xvec$ witnesses that $f_{\<\evec,\kkvec\>}$ is active.

 Let $K'\in \br K;{\omega};$,
$\adot=\{\<p_{\xi},x_{\xi}\>:{\xi}\in K'\}$ and 
 ${\zeta}\in K\setm K'$.
Let $r^*=\<c^*,q^*\>\le p_{\zeta}$ be a determined condition such that 
for each active $f_{\<\evec,\kkvec\>}$ and $w\in Z$
the value
$f^{q^*}_{\<\evec,\kkvec\>}(w)$ is defined.
Let
\begin{displaymath}
Y=\{f^{r^*}_{\<\evec,\kkvec\>}(w):\text{$f_{\<\evec,\kkvec\>}$ is active 
and $w\in Z$}\}.
\end{displaymath}

\begin{pClaim}
$Y$ is finite.
\end{pClaim}

\begin{proof}[Proof of the claim]
Since $\{f_{\eta}:{\eta}\in J\}$ acts loop-free on $Z_{\zeta}\setm Z$,
the elements of a witnessing sequence are pairwise different, so there
are only finitely many of them and a witnessing sequence works only
for one active  $f_{\<\evec,\kkvec\>}$. So there is only finitely
many active  $f_{\<\evec,\kkvec\>}$. 
\end{proof}

We show
that 
\begin{equation}
\tag{$\bullet$}
r^*\force \adot\text{ is dense in $G[\oo\setm Y]$.}
\end{equation}
which completes the proof of lemma \ref{lm:many-dense}.

To verify $(\bullet)$ assume that $r'\le r^*$, $r'=\<c',q'\>$ is determined,
$B\in \br \oo\setm Y;<{\omega};$ and $b\in {}^B2$.

Pick ${\xi}\in K$ such that $\supp(c')\cap \supp (c_{\xi})=Z$
and $\dom(q')\cap \dom(q_{\xi})=J$. 
To prove $(\bullet)$ 
it is enough 
to construct a common extension $p=\<c,q\>$ of
$r'=\<c',q'\>$ and $p_{\xi} =\<c_{\xi} ,q_{\xi} \>$
such that $c(x_{\xi},{\beta})=b({\beta})$ for each ${\beta}\in B$.

Let $\supp c =\supp c' \cup\supp c_{\xi} $.
Put $\dom q=\dom q'\cup\dom q_{\xi}$ and let
\begin{displaymath}
q({\eta})= \left\{ 
\begin{array}{ll}
q'({\eta})\cup q_{\xi}({\eta})&\mbox{if ${\eta}\in J$},\\
q'({\eta} )&\mbox{if ${\eta}\in \dom q'\setm J$,}\\
q_{\xi}({\eta})&\mbox{if ${\eta}\in \dom q_{\xi}\setm J$.}
\end{array}
\right.
\end{displaymath}
Put $c^-=c'\cup c_{\xi}$.

We should define $c\supset c^-$ on the set 
$$E=\{\{a,b\}:a\in Y_{\xi}\setm Y,b\in \supp c'\setm Y\}.
$$ 
such that every $q({\eta})$ is a partial isomorphism of $G$,
more precisely, $q({\eta})\in \isop(G^c,G^c)$.
To do so, observe that if we take 

$$E^+=\{\{a,b\}:a\in Y_{\xi}\setm Y,b\in \supp c'\}
$$ 
and
for $e\in E^+$ define $a_e=e\cap (Z_{\xi}\setm Z)$ and 
$b_e=e\cap \supp c'$ then
$q({\eta})\in \isop(G^c)$  if and only if $(\dag)$ below holds:
\begin{equation}\tag{\dag}
\text{if $e=\{a_e,b_e\}\in E^+$ then $c\{a_e,b_e\}=
c\{q_{\xi}({\eta})(a_e),q'({\eta})(b_e)\}$}.
\end{equation}

Define an equivalence relation  ${\equiv}$ on $E^+$:
$e {\equiv}  e'$  iff $e=e'$ or there is an active 
$f_{\<\evec,\kkvec\>}$ such that 
  $a_{e'}=f^{q_{\xi}}_{\<\evec,\kkvec\>}(a_{e})$ 
and $b_{e'}=f^{q'}_{\<\evec,\kkvec\>}(b_{e})$.

\begin{claim}\label{cl:1}
If $e\equiv e'$ and $a_e=a_{e'}$ then $e=e'$.
\end{claim}

\begin{proof}[Proof of the claim \protect\ref{cl:1}]
Assume $e {\equiv} e'$ and $b_e\ne b_{e'}$.
Then there is an active  
$f_{\<\evec,\kkvec\>}$ such that 
  $a_{e'}=f^{q_{\xi}}_{\<\evec,\kkvec\>}(a_{e})$ 
and $b_{e'}=f^{q'}_{\<\evec,\kkvec\>}(b_{e})$.
Since
$1\force$ ``$\{f_{\eta}:{\eta}\in J\}$ acts freely on $\oo\setm {\mu}$''
it follows that  $a_{e}\ne f^{q'}_{\<\evec,\kkvec\>}(a_{e})$ and
so $a_e\ne a_{e'}$.
\end{proof}

\begin{claim}\label{cl:2}
If $e,e'\in E^+\cap \dom (c^-)$ and $e\equiv e'$ then 
$c^-(e)=c^-(e')$. 
\end{claim}

\begin{proof}[Proof of the claim \protect\ref{cl:2}]
Fix an active 
$f_{\<\evec,\kkvec\>}$ such that 
  $a_{e'}=f^{q_{\xi}}_{\<\evec,\kkvec\>}(a_{e})$ 
and $b_{e'}=f^{q'}_{\<\evec,\kkvec\>}(b_{e})$.
Since $e,e'\in E^+\cap \dom (c^-)$ it follows that 
$e,e'\in\dom (c_{\xi})$ and so $a_e,a_{e'}\in Z$.
Thus $b_{e'}=f^{q_{\xi}}_{\<\evec,\kkvec\>}(b_{e})$.
But $f^{q_{\xi}}_{\<\evec,\kkvec\>}\in\isop(G^{c_{\xi}},G^{c_{\xi}})$ for
$\<c_{\xi},q_{\xi}\>\in\ccal*P_{\alpha}$ so 
$c_{\xi}(e)=c_{\xi}(e')$.
\end{proof}

\begin{claim}\label{cl:3}
If $e\in E^+\cap \dom(c^-)$ and $e\equiv e'$ then
$b_{e'}\in Y$. 
\end{claim}
\begin{proof}[Proof of the claim \protect\ref{cl:3}]
Since $e\in E^+\cap \dom(c^-)$ we have $b_e\in Z$.
Fix an active  
$f_{\<\evec,\kkvec\>}$ such that 
  $a_{e'}=f^{q_{\xi}}_{\<\evec,\kkvec\>}(a_{e})$ 
and $b_{e'}=f^{q'}_{\<\evec,\kkvec\>}(b_{e})$.
Since $f_{\<\evec,\kkvec\>}$  is active 
it follows that $f^{q^*}_{\<\evec,\kkvec\>}(b_e)$ is
defined and $f^{q^*}_{\<\evec,\kkvec\>}(b_e)\in Y$.
But
$f^{q'}_{\<\evec,\kkvec\>}(b_e)=
f^{q^*}_{\<\evec,\kkvec\>}(b_e)$
so $b_{e'}\in Y$ which was to be proved.
\end{proof}

By claims \ref{cl:1}--\ref{cl:3} 
 we can find a condition $c\in\ccal$ with 
$\supp c= \supp c'\cup\supp c_{\xi}$ and $\dom c =\br \supp c ; 2;$
such that
\begin{enumerate}\alabel
\item $c\supset c^-=c'\cup c_{\xi}$,
\item $c(e)=c(e')$ whenever $e {\equiv} e'$,
\item $c\{x_{\xi},{\beta}\}=b({\beta})$ for ${\beta}\in B$. 
\end{enumerate}
Then $(\dag)$ holds and as we have seen above, $\<c,q\>\in
\ccal*P_{\alpha}$ and 
\begin{displaymath}
\<c,q\>\force (\forall {\beta}\in B)\ \{x_{\xi},{\beta}\}\in E(G)
\text{ iff }b({\beta})=1.
\end{displaymath}
Thus $(\bullet)$ holds. Hence lemma \ref{lm:many-dense} is proved.
\end{proof}  
So we have shown that (II) is preserved during the inductive construction,
which was the last step to prove theorem \ref{tm:smooth}.
\end{proof}

\end{document}